\documentclass[a4paper,10pt]{article}

\usepackage[english]{babel}
\usepackage{amssymb,amsmath}

\newtheorem{thm}{Theorem}[section]
\newtheorem{prop}[thm]{Proposition}
\newtheorem{defi}[thm]{Definition}
\newtheorem{lem}[thm]{Lemma}

\def\Qset{\mathbf{Q}}
\def\Fset{\mathbf{F}}

\def\Supp{\mathrm{Supp}}

\def\pf{{\it Proof. }}
\def\rmk{{\it Remark. }}

\newcommand{\St}[3]{#1\mathrm{#2}#3}

\title{Solving Fermat-type equations $x^5+y^5=dz^p$}
\author{Nicolas Billerey\thanks{billerey@math.jussieu.fr}\,
 and\, Luis V. Dieulefait\thanks{ldieulefait@ub.edu}}

\begin{document}

\maketitle

\begin{abstract}
In this paper, we are interested in solving the Fermat-type equations $x^5+y^5=dz^p$ where $d$ is a positive integer and $p$ a prime number $\ge 7$. We describe a new method based on modularity theorems which allows us to improve all the results of~\cite{Bil07}. We finally discuss the present limitations of the method by looking at the case $d=3$.   
\end{abstract}

		\section{Introduction}

Let $p$ be a prime number $\ge 7$ and $d$ be a positive integer. We say that
a solution $(a,b,c)$ of the equation $x^5+y^5=dz^p$ is primitive if $(a,b)=1$
and non-trivial  if $c\not=0$ (note that this is not the same
definition as in \cite{Bil07}). Let us recall briefly the
generalization of the so-called modular method of Frey for solving this equation.

Assume that $(a,b,c)$ is a non-trivial  primitive solution of $x^5+y^5=dz^p$. Then the equation
\renewcommand{\theequation}{$\star$}
\begin{equation}\label{eq:FreyCurve}
 y^2=x^3-5(a^2+b^2)x^2+5\left(\frac{a^5+b^5}{a+b}\right)x
\end{equation}
defines an elliptic curve $E(a,b)$ over $\Qset$ of conductor $N$ which is semistable at each prime different from $2$ and $5$.
By results of Wiles, Taylor-Wiles, Diamond and Skinner-Wiles, $E(a,b)$ is modular.
Furthermore, $E(a,b)$ is a Frey-Hellegouarch curve in the following sense~: the Galois
representation $\rho_p$ on $p$-torsion points of $E(a,b)$ is irreducible and unramified outside
$2$, $5$, $p$ and the set of primes dividing $d$. The conductor $N(\rho_p)$ (prime to
$p$) and the weight $k$ of $\rho_p$ are computed in \cite[\S3]{Bil07}. Thus, it follows from a theorem of Ribet that there exists a
modular form $f$ of weight $k$, level $N(\rho_p)$ and trivial character such that the
associated $p$-adic representation $\sigma_{f,p}$ satisfy
$\sigma_{f,p}\equiv \rho_p\pmod{p}$. More precisely, let us denote  by
$a_q$ and $a_q'$ the coefficients of the $L$-functions of $E$ and $f$ respectively,
by $K_f$ the number field generated by all the $a_q'$'s numbers and by $N^{K_f}_{\Qset}$ the corresponding
norm map. We then have the following proposition.
\begin{prop}
There exists a primitive newform $f$ of weight $k$ and level $N(\rho_p)$ such that, for
each prime $q$, the following conditions hold.
\begin{enumerate}\label{prop:Ribet}
    \item If $q$ divides $N$ but $q$ does not divide $pN(\rho_p)$, then
\begin{equation*}
p \textrm{ divides } N^{K_f}_{\Qset}\left(a_q'\pm(q+1)\right).
\end{equation*}
    \item If $q$ does not   divide $pN$, then
\begin{equation*}
p \textrm{ divides } N^{K_f}_{\Qset}\left(a_q'-a_q\right).
\end{equation*}
\end{enumerate}
\end{prop}
The aim of the modular method is to contradict the existence of such a form~$f$. We
describe, in the following section, a method which allows us sometimes to reach
this goal.

\addtocounter{equation}{-1}
\renewcommand{\theequation}{\arabic{equation}}

        \section{Description of the method}\label{s:method}

Assume that $(a,b,c)$ is a non-trivial  primitive solution of $x^5+y^5=dz^p$. Let $f$ be the newform of Prop.~\ref{prop:Ribet} and $q$ be a prime number which does not
divide $pN(\rho_p)$. Assume moreover that
$p$ does not divide $N^{K_f}_{\Qset}\left(a_q'\pm(q+1)\right)$. Then $E(a,b)$ has good
reduction at~$q$. So, we are reduced to contradict the second assertion of the proposition.
Since $E(a,b)$ has a non-trivial  $2$-torsion group, the coefficient $a_q$ is even.
Furthermore, since the equation~(\ref{eq:FreyCurve}) of $E(a,b)$ does not depend on $p$, we can look at the reduction modulo~$q$ of the Frey curve without caring
whether $(a,b,c)$ is a solution of the Fermat equation or not. A short calculation (for instance, on \verb"pari/gp") gives us the list of possibilities
 for $a_q$,
 when $(a,b)$ describes $\Fset_q\times \Fset_q$. Assume now that $p$ does not   divide
 $N^{K_f}_{\Qset}\left(a_q'-a_q\right)$, for any $a_q$ in this list. Contradiction
 will follow if, for each form $f$ as above, we can find such a prime~$q$.

        \section{Applications to the Fermat equation}

We apply, in this section, the method described above to some values of $d$.

	\subsection{Case where $d=2^{\alpha}\cdot 3^{\beta}\cdot 5^{\gamma}$}

In this paragraph, we are interested in the case where
\begin{equation*}
 d=2^{\alpha}\cdot 3^{\beta}\cdot 5^{\gamma},\quad \text{with } \alpha\ge 2 \text{ and }
 \beta,  \gamma \text{ arbitrary}.
\end{equation*}
The following theorem generalizes Theorems~$1.2$ and $1.3$ of~\cite{Bil07}.
\begin{thm}
Assume $d$ is as above. Then, the equation $x^5+y^5=dz^p$ does not have non-trivial
primitive solutions for $p\ge 13$.
\end{thm}
\pf Assume that $(a,b,c)$ is a non-trivial  primitive solution. It follows
from~\cite[\S$3$]{Bil07}, that the representation $\rho_p$ is irreducible of weight
$k=2$.

If $\beta=0$, then we have $N(\rho_p)=25$ or $50$. Since there is no newform of weight $2$
and level $25$, we necessary have $N(\rho_p)=50$. There are exactly two such forms and both
of them have rational coefficients. The curve $E(a,b)$ is semistable at $q=3$. Assume that
$E(a,b)$ has multiplicative reduction at $3$. By Prop.~\ref{prop:Ribet}, we have
$a_3'\pm 4\equiv 0\pmod{p}$. Besides, by~\cite{mfd}, we have $a_3'=\pm 1$ which is a
contradiction, since $p\ge 13$. So $E(a,b)$ has good reduction at $q=3$ and by the
proposition above, $\pm 1=a_3'\equiv a_3\pmod{p}$. This is also a contradiction because
$a_3$ is even ($E(a,b)$ has a non-trivial $2$-torsion subgroup) and
$\lvert a_3\rvert\le 2\sqrt{3}$, \emph{i.e.} $a_3=0$ or $\pm2$.

If $\beta > 0$, then we have $N(\rho_p)=75$ or $150$. Assume that we have $N(\rho_p)=75$. By~\cite{mfd}, there are
exactly $3$ primitive
newforms of weight $2$ and level $75$. They all have coefficients in~$\Qset$ and the
form $f$ of Prop.~\ref{prop:Ribet} is one of them. Moreover, by~\cite{mfd},
 we have $a_7'=0$ or $\pm 3$. Since $p\ge 13$, the first condition of Prop.~\ref{prop:Ribet}
  does not hold for $q=7$ and $E(a,b)$ has good reduction at $7$. Following the method described
  in the previous section, we find that $a_7$ belongs to the set $\left\{-4,-2,2\right\}$.
  We then deduce that the second condition of Prop.~\ref{prop:Ribet} does not hold too. In
  other words, we have $N(\rho_p)=150$.

There are exactly $3$ primitive newforms of weight $2$ and level $150$,
denoted by $\St{150}{A}{1}$, $\St{150}{B}{1}$ and $\St{150}{C}{1}$ and $f$ is one of
them. If $f=\St{150}{B}{1}$, then $a_7'=4$ and a contradiction follows as above. So, $f=\St{150}{A}{1}$ or $\St{150}{C}{1}$ and by~\cite{mfd}, we have $a_{11}'=2$.
 Since $p\ge 13$, the first condition of Prop.~\ref{prop:Ribet} does not hold for $q=11$ and $E(a,b)$
 has good reduction at $11$. Besides, we have $a_{11}=0$ or $\pm 4$. So, the second
 condition of Prop.~\ref{prop:Ribet} does not hold too and we obtain a contradiction. This
 ends the  proof of the theorem.

	\subsection{Case where $d=7$}

In this paragraph, we prove the following theorem.
\begin{thm}
The equation $x^5+y^5=7z^p$ does not have non-trivial
primitive solutions for $p\ge 13$.
\end{thm}
\pf Assume that $(a,b,c)$ is a non-trivial  primitive solution. It follows
from~\cite[\S$3$]{Bil07}, that the representation $\rho_p$ is irreducible, of weight
$k=2$ (since $p\not=7$) and level $N(\rho_p)=350$, $1400$ or $2800$.

Let us first assume that the form $f$ of Prop.~\ref{prop:Ribet} has eigenvalues
which are not rational integers. There are exactly $19$ such forms and for all of them
we have $a_3'=\alpha$ where $\alpha$ is the generator
 of the field $K_f$ given in~~\cite{mfd}. If $E(a,b)$ has good reduction at $q=3$,
 we have $a_3=\pm 2$. Furthermore,  $N^{K_f}_{\Qset}
 \left(a_3'\pm 2\right)$ belong to the set $\left\{\pm 2,\pm 4, -6,\pm 10\right\}$.
Since $f$
 satisfies the second condition of Prop.~\ref{prop:Ribet}, we deduce that $E(a,b)$ has
 multiplicative reduction at $3$.

If $f$ is not one the forms denoted by $\St{1400}{S}{1}$, $\St{1400}{T}{1}$,
$\St{2800}{QQ}{1}$ or $\St{2800}{RR}{1}$ in~\cite{mfd}, then $N^{K_f}_{\Qset}
\left(a_3'\pm 4\right)$ belong to $\left\{4,8,10,12,16,20\right\}$ and $p$ divides one of them. This is a
contradiction. So, $f$ is necessary one of the $4$ forms above and we have $N^{K_f}_{\Qset}
\left(a_3'\pm 4\right)=\pm 2\cdot 29$ or $\pm 2\cdot 11$. It then follows that $p=29$.
Besides, if $E(a,b)$ has good reduction at $q=17$, then
$a_{17}\in\left\{0,2,4,\pm 6,-8\right\}$, but by~\cite{mfd}, $29$ does not divide
$N^{K_f}_{\Qset}\left(a_{17}'\right)$, $N^{K_f}_{\Qset}\left(a_{17}'-2\right)$,
$N^{K_f}_{\Qset}\left(a_{17}'-4\right)$, $N^{K_f}_{\Qset}\left(a_{17}'\pm 6\right)$
and $N^{K_f}_{\Qset}\left(a_{17}'+8\right)$. So,  $E(a,b)$
has multiplicative reduction at $q=17$ and $29$ divides $N^{K_f}_{\Qset}\left(a_{17}'
\pm 18\right)=\pm 2^6\cdot 79$ or $\pm 2^4\cdot 359$. This leads us again to a contradiction and we
conclude that the eigenvalues of $f$ are all rational integers.

In other words, $f$ corresponds to an elliptic curve defined over~$\Qset$. There are
exactly $6$ isogeny classes of elliptic curves of level $350$, $14$ of level $1400$
and $33$ of level $2800$. For all of them, we will contradict the conditions of
Prop.~\ref{prop:Ribet} with $q=3$, $11$, $19$, $23$ or $37$. As we have seen
in~\S\ref{s:method}, if $E(a,b)$ has good reduction at $q$, we can list the possible
values of $a_q$. For the prime numbers $q$ above, we find
\begin{equation*}
a_3=\pm 2, \qquad a_{11}\in\lbrace0,\pm 4\rbrace, \qquad a_{19}\in\lbrace0,\pm 4\rbrace,
\end{equation*}
\begin{equation*}
a_{23}\in\lbrace 0,\pm 2, \pm 4, \pm 6, \pm 8\rbrace\quad\text{and}\quad  a_{37}\in
\lbrace 0,-2,\pm 4, -6, \pm 8,\pm 10, 12\rbrace.
\end{equation*}
By the Hasse-Weil bound, $E(a,b)$ has good reduction at $q=3$. We then deduce that
$f$ satisfies $a_3'=\pm 2$. Among these curves, let us begin to deal with those without
$2$-torsion rational over $\Qset$. If $f$ is one of the curves denoted by $\St{2800}{W}{1}$ and
$\St{2800}{AA}{1}$ in~\cite{mfd}, we have $a_{11}'=\pm3$ and this contradicts the
congruences of Prop.~\ref{prop:Ribet} with $q=11$. If $f$ is one of the curves denoted by
$\St{1400}{D}{1}$, $\St{1400}{K}{1}$, $\St{2800}{D}{1}$ and
$\St{2800}{N}{1}$, we have $a_{11}'=\pm 1$. We then have a contradiction except maybe for
 $p=13$. Besides, for these $4$ curves, we have $a_{23}'=\pm 3$
and the same argument implies another contradiction except for $p=19$. Brought together, these two
results imply that $f$ is not one of these $4$ forms. If now $f$ is one of the curves denoted by
$\St{1400}{C}{1}$, $\St{1400}{N}{1}$, $\St{2800}{E}{1}$ and
$\St{2800}{M}{1}$, we have $a_{11}'=\pm 5$. We then have a contradiction except maybe for
$p=17$. Besides, for these curves, we have $a_{19}'=\pm 2$. By the same argument as before, it then follows
once more a contradiction.

The two remaining curves of level $350$, $1400$ or $2800$ such that $a_3'=\pm 2$, denoted by
$\St{1400}{H}{1}$ and $\St{2800}{G}{1}$ are the
only two curves, with non-trivial $2$-torsion group. They satisfy $a_{19}'=\pm 2$ and
$a_{37}'=6$. Since these values do not belong to the set of possible values for $a_{19}$
and $a_{37}$ described above, we finally have a contradiction to the existence of a
non-trivial primitive solution of $x^5+y^5=7z^p$.

	\subsection{Case where $d=13$}

In this paragraph, we prove the following theorem.
\begin{thm}
The equation $x^5+y^5=13z^p$ does not have non-trivial
primitive solutions for $p\ge 19$.
\end{thm}
\pf Assume that $(a,b,c)$ is a non-trivial  primitive solution. It follows
from~\cite[\S$3$]{Bil07}, that the representation $\rho_p$ is irreducible, of weight
$k=2$ (since $p\not=13$) and level $N(\rho_p)=650$, $2600$ or $5200$.

Let $q$ be a prime number different from $2$, $5$, $13$ and $p$. By Prop.~\ref{prop:Ribet}, $p$ divides
either $N^{K_f}_{\Qset}\left(a_q'\pm(q+1)\right)$ or $N^{K_f}_{\Qset}\left(a_q'-a_q\right)$. In other words, $p$
is a prime factor of the resultant $R_q$ of the minimal polynomial of $a_q'$ and $P_q(X)=(X^2-(q+1)^2)\prod(X-a_q)$ where
the product runs over all possible values of $a_q$. For instance, if $q=3$, then $P_3(X)=(X^2-16)(X^2-4)$.

Let us first assume that $f$ has rational Fourier coefficients. If $a_3'\not=\pm2$, then $R_3$ has only $2$, $3$, $5$ and $7$ as prime factors.
So, we deduce that $a_3'=\pm 2$. There are exactly $6$ such newforms of level $650$, $5$ of level $2600$ and $37$ of level 
$5200$ (for the curves of level $5200$, the notation will exceptionally refer to~\cite{Cr97}) . For all of
them, $a_{7}'$ does not belong to the list $\lbrace\pm 2,-4\rbrace$ of possible values for $a_7$ when $E(a,b)$ has good
reduction at $7$. The same observation
holds for the $13$ elliptic curves of level $5200$ with $a  _3'=\pm 2$ except for those denoted by $\St{5200}{S}{1}$,
$\St{5200}{BB}{1}$, $\St{5200}{AA}{1}$ and $\St{5200}{Z}{1}$ (in~\cite{Cr97}). If $f$ is one of the first three of them, then we have $a_{11}'=6$ or
$\pm2$. Besides, if $E(a,b)$ has good reduction at $11$, then $a_{11}$ belongs to $\lbrace0,\pm 4\rbrace$. So, this is
a contradiction and $f=\St{5200}{Z}{1}$. Nevertheless, in this case, $a_{17}'=-2$ does
not belong to the set $\lbrace0,2,4,\pm6,-8 \rbrace$ of possible values for $a_{17}$ when $E(a,b)$ has good reduction at $17$.
We then deduce that the Fourier coefficients of $f$ are not all rational.

Let us now assume that $N(\rho_p)=650$ or $2800$. For each $f$ in these levels, $a_3'=\alpha$ is a root of the polynomial
defining $K_f$ given in \cite{mfd}. We then verify that $R_3$ is supported only by $2$ and $5$ except for the curves 
denoted $\St{2800}{QQ}{1}$ and $\St{2800}{RR}{1}$. But, they both satisfy $a_7'=\pm 1$ and this leads us to a 
contradiction.

So, we necessarily  have $N(\rho_p)=5200$. There are exactly $29$ newforms of this level with non-rational eigenvalues numbered 
from 38 to 66. Four of them (those numbered $39$, $42$, $46$ and $47$) satisfy $a_3'=0$ or $\pm 1$. 
So, $f$ is not one of them. If $f$ is the curve numbered $63$, then 
the field of coefficients is generated by a root $\alpha$ of the following polynomial 
$x^4 + 6x^3 - 18x^2 - 30x + 25$ and 
\begin{equation*}
a_3'=\frac{1}{10}\left(\alpha^3+6\alpha^2-13\alpha-20\right).
\end{equation*}
Its characteristic polynomial is then $x^4+2x^3-7x^2-8x+16$ and we get $R_3=2^{18}$ in this case. This is of course a 
contradiction. The same conclusion will follow if $f$ is the curve numbered $64$, since, in this case, the generating polynomial is $x^4 + 6x^3 - 87x^2 - 492x + 604$ and the characteristic polynomial of $a_3'$ is 
 $x^4-2x^3-7x^2+8x+16$.

For all the other curves, $a_3'=\alpha$ is a root of the generating polynomial  of $K_f$ given in the tables and we have a contradiction in the same way as before by looking at $R_3$ except 
for the following eight pairs $(f,p)$~:
\begin{equation*}
(f=54,p=43),\quad (f=55,p=43),\quad (f=58,p=23),\quad (f=59,p=67),  
\end{equation*}
\begin{equation*}
\quad (f=61,p=23),\quad (f=62,p=67),\quad (f=65,p=23),\quad (f=66,p=43).
\end{equation*}
For all of them, we have a contradiction as before by looking at the coefficient $a_7'$ except for the last two curves 
where we have to consider $a_{19}'$. 

We finally deduce a contradiction to the existence of a non-trivial primitive solution of the equation $x^5+y^5=13z^p$.

        \section{The case $d=3$ and limitations of the method}

It is clear that the method will not work if there exists an el\-liptic
curve over $\Qset$ of the form~(\ref{eq:FreyCurve}) and level $N(\rho_p)$ (for large $p$). For
convenience, we adopt the fol\-lowing definition which makes this observation precise (where $\Supp$ denotes the support of an integer and $v_2$ the $2$-adic valuation of~$\Qset$).
\begin{defi}\label{defi:modobstruction}
We say that there is a modular obstruction for the equation $x^5+y^5=dz^p$ (or just
for $d$) if there exists $(a,b)$ two coprime integers such that the following two
conditions hold.
\begin{enumerate}
    \item The integer $m=a^5+b^5$ is non-zero and we have
\begin{equation*}
\Supp(m)\setminus\left\{2,5\right\}=\Supp(d)\setminus\left\{2,5\right\}.
\end{equation*}
    \item We have~:
\begin{itemize}
        \item if $\Supp(d)$ is not included in $\lbrace 2,5\rbrace$, then $ab\not=0$,
	\item if $\Supp(d)$ is included in $\lbrace 2,5\rbrace$ and $d$ is even, then $ab\not=0$,
        \item if $d$ is odd, then $v_2(m)\not= 2$,
        \item if $v_2(d)=1$, then either $v_2(m)\ge 3$, or $v_2(m)=1$, or $v_2(m)=0$ and
        $\max(v_2(a),v_2(b))=1$,
        \item if $v_2(d)=2$, then $v_2(m)=2$,
        \item if $v_2(d)\ge 3$, then $v_2(m)\ge 3$.
\end{itemize}
\end{enumerate}
\end{defi}

The following lemma gives a sufficient condition which insures that
there is no modular obstruction, for several $d$.
\begin{lem}
Let $d$ be a positive integer such that for any prime $\ell$ dividing $d$, we have
$\ell\not\equiv 1\pmod{5}$. Then, there is a modular obstruction for $d$ if and only if
$d=5^{\gamma}$ or $d=2\cdot 5^{\gamma}$ with $\gamma \ge 0$.
\end{lem}
\pf Assume that there is a modular obstruction for $d$ given by two coprime integers $(a,b)$.
Then $m=a^5+b^5$ is non zero and $\Supp(m)\setminus\left\{2,5\right\}=\Supp(d)\setminus
\left\{2,5\right\}$. Following \cite{Bil07}, let us denote by $\phi$ the irreducible polynomial
\begin{equation*}
\phi(x,y)=x^4-x^3y+x^2y^2-xy^3+y^4.
\end{equation*}
By Lemmas~$2.5$ and $2.6$ of~\cite{Bil07} and the hypothesis, we have~:
\begin{enumerate}
	\item either $5$ divides $m$ and then $5$ divides $a+b$ and $\phi(a,b)=\pm 5$;
	\item or $5$ does not divide $m$ and then $\phi(a,b)=\pm 1$.
\end{enumerate}
In other words, $(a,b)$ is a solution of a Thue equation of the form $\phi(x,y)=A$, where $A=\pm1$ 
or $\pm 5$ and we can assume that $a\not=0$ ($\phi$ is symmetric). Since $\phi$ is totally complex, 
 this leads to 
\begin{equation*}
\lvert A\rvert = \lvert a\rvert^4\prod_{k=1}^{4}\lvert b/a-\alpha_k\rvert\ge \lvert a\rvert^4 
\sin^2\left(\frac{2\pi}{5}\right)\cdot\sin^2\left(\frac{4\pi}{5}\right)\ge 0.312\cdot \lvert a\rvert^4,
\end{equation*}
where $\alpha_k=-\exp(2ik\pi/5)$, $k=1,\ldots,4$, are the roots of $\phi$. This gives an upper bound 
for $\lvert a\rvert$.

In the
first case, it implies that we have $(a,b)=(1,-1)$ or $(-1,1)$ and then $m=0$, which is a contradiction.
In the second case, we deduce
\begin{equation*}
(a,b)\in\lbrace(1,1),(-1,-1),(\pm 1,0),(0,\pm 1)\rbrace.
\end{equation*}
In other words, $m=\pm 1$ or $m=\pm 2$. By the first condition of
Def.~\ref{defi:modobstruction}, there exists $\alpha,\gamma \ge 0$ such that
$d=2^{\alpha}\cdot 5^{\gamma}$. Since $v_2(m)=0$ or $1$, we have, by the second condition,
$\alpha=0$ or $1$.

Conversely, if $d=5^{\gamma}$ or $d=2\cdot 5^{\gamma}$ with $\gamma \ge 0$, there is a
modular obstruction for $d$ given, for example, by $(a,b)=(1,1)$.

\medskip

\rmk For $d=11$, there is a
modular obstruction given by the elliptic curves $E(2,3)$ or $E(3,-1)$. Note that finding a modular obstruction for a given $d$ involves solving some Thue-Mahler equation which, at least theoretically, is possible but can be difficult in practice.

\medskip

Let us now look at the case where $d=3$. By the previous lemma, there is no
modular obstruction. Nevertheless, as we will see, we were not able to solve this equation for all $p$.

Fix for now a prime $p$ and $(a,b,c)$ a non-trivial primitive solution of the equation $x^5+y^5=3z^p$. 
The following lemma makes more precise Lemma~$4.3$ of \cite{Bil07}. We warn the reader that in this paragraph we are using 
only Stein's notations \cite{mfd} for modular forms (including elliptic curves) which is not the case in \cite{Bil07} 
where the author was referring to Cremona's Tables of elliptic curves \cite{Cr97}.
\begin{lem}\label{lem:remainingcurves}
If $p\ge 17$, then we have
\begin{enumerate}
    \item either $5$ divides $a+b$ and $f=\St{1200}{K}{1}$,
    \item or $5$ does not divide $a+b$ and $f=\St{1200}{A}{1}$.
\end{enumerate}
\end{lem}
\pf Assume that $5$ divides $a+b$. By Lemma~$4.3$ of \cite{Bil07}, $f$ is one of following newforms (with Stein's notations)~:
\begin{equation*}
 \St{150}{B}{1}, \St{600}{C}{1}, \St{600}{A}{1}, \St{1200}{O}{1}, \St{1200}{R}{1}, \St{1200}{E}{1}, \St{1200}{K}{1}.
\end{equation*}
If $f=\St{150}{B}{1}$, $\St{600}{C}{1}$, $\St{1200}{O}{1}$, $\St{1200}{R}{1}$ or $\St{1200}{E}{1}$, we have $a_7'=0$
or $4$. Besides, if $E(a,b)$ has good reduction at $7$, we have $a_7=\pm2$ or $-4$. We then obtain a contradiction by
looking at the  conditions of Prop.~\ref{prop:Ribet} for $q=7$. If $f=\St{600}{A}{1}$, then $a_{13}'=6$. Besides, if $E(a,b)$ has good reduction at $13$, then $a_{13}$ belongs to the set $\lbrace0,\pm 2,\pm 4\rbrace$. So, there is again a contradiction. So, $f=\St{1200}{K}{1}$ in this case.

Assume now that $5$ does not divide $a+b$. By Lemma~$4.3$ of \cite{Bil07}, $f$ is one of following newforms
(with Stein's notations)~:
\begin{equation*}
 \St{150}{A}{1}, \St{150}{C}{1}, \St{600}{D}{1}, \St{600}{G}{1}, \St{1200}{H}{1}, \St{1200}{L}{1}, \St{1200}{G}{1},
 \St{1200}{A}{1}, \St{1200}{M}{1}, \St{1200}{S}{1}.
\end{equation*}
For $f=\St{1200}{S}{1}$ we have $a_7'=4$ and using this coefficient we derive again a contradiction. For all the other curves
except $\St{1200}{A}{1}$, we have $a_{11}'=\pm 2$. Besides, if $E(a,b)$ has good reduction at $11$, we have $a_{11}=0$
or $\pm 4$. So, $f$ is not one of them and we conclude that $f=\St{1200}{A}{1}$ in this case.

\medskip

If $f=\St{1200}{K}{1}$ or $\St{1200}{A}{1}$, then for any prime $q>5$ smaller than $5000$, the Fourier coefficient $a_q'$
of $f$ actually lies in the list of possible values for $a_q$. This is why we have not been able to prove the emptiness of
the set of non-trivial primitive solutions for $d=3$.

Nevertheless, we will give a criterion which allows us to conclude for a fixed $p$ and verify that it holds for any
$17\le p\le 10^6$. Let us consider $q$ a prime number congruent to $1$ modulo $p$ and write $q=np+1$. The group $\mu_n(\Fset_q)$ of $n$th roots of unity in $\Fset_q$ has order $n$. We now define four subsets $A^{\pm}(n,q)$ and $B^{\pm}(n,q)$ of $\mu_n(\Fset_q)$ in the following way.

\begin{enumerate}
\item Let $\widetilde{A}(n,q)$ be the subset of
$\mu_n(\Fset_q)$ consisting of all $\zeta$ such that
\begin{equation*}
       405+62500\zeta\ \textrm{is a square in}\ \Fset_q.
\end{equation*}
For such a $\zeta$, let us consider the smallest integer $\delta_{1,\zeta}\ge 0$ such that
\begin{equation*}
\delta_{1,\zeta}^2\pmod{q}= 405+62500\zeta. \\
\end{equation*}
We define $A^{+}(n,q)$ (resp. $A^{-}(n,q)$) as the subset of $\widetilde{A}(n,q)$ consisting of $\zeta$ such that
\begin{equation*}
	-225+10\delta_{1,\zeta}\quad(\textrm{resp. } -225-10\delta_{1,\zeta})
\end{equation*}
is a square modulo $q$. For any $\zeta\in A^+(n,q)$, let us consider the cubic curve over $\Fset_q$ defined by the 
following equation
\begin{equation*}
F_{1,\zeta}^+: y^2=x^3-\frac{\delta_{1,\zeta}}{25}x^2+25\zeta x.
\end{equation*}
Its discriminant $6480\zeta^2=2^4\cdot3^4\cdot5\zeta^2$ is non-zero and $F_{1,\zeta}^+$ is an elliptic curve over $\Fset_q$. Let us denote by $n_{1,q}^+(\zeta)$ the number of $\Fset_q$-rational points of $F_{1,\zeta}^+$ and write
\begin{equation*}
    a_q^+(\zeta)=q+1-n_{1,q}^+(\zeta).
\end{equation*}
If $\zeta\in A^{-}(n,q)$, let us define in the same way, the cubic curve
\begin{equation*}
F_{1,\zeta}^-: y^2=x^3+\frac{\delta_{1,\zeta}}{25}x^2+25\zeta x.
\end{equation*}
As a twist of $F_{1,\zeta}^+$, it is also an elliptic curve over $\Fset_q$ and we write
\begin{equation*}
    a_q^-(\zeta)=q+1-n_{1,q}^-(\zeta),
\end{equation*}
where $n_{1,q}^-(\zeta)$ denotes the number of $\Fset_q$-rational points of $F_{1,\zeta}^-$.

\item Let $\widetilde{B}(n,q)$ be the subset of
$\mu_n(\Fset_q)$ consisting of all $\zeta$ such that
\begin{equation*}
       405+20\zeta\ \textrm{is a square in}\ \Fset_q.
\end{equation*}
For such a $\zeta$, let us consider the smallest integer $\delta_{2,\zeta}\ge 0$ such that
\begin{equation*}
\delta_{2,\zeta}^2\pmod{q}= 405+20\zeta. \\
\end{equation*}
We define $B^+(n,q)$ (resp. $B^-(n,q)$) as the subset of $\widetilde{B}(n,q)$ consisting of $\zeta$ such that
\begin{equation*}
	-225+10\delta_{2,\zeta}\quad(\textrm{resp. } -225-10\delta_{2,\zeta})
\end{equation*}
is a square modulo $q$. For any $\zeta\in B^+(n,q)$, let us consider the cubic curve over $\Fset_q$ defined by the following equation
\begin{equation*}
F_{2,\zeta}^+: y^2=x^3-\delta_{2,\zeta}x^2+5\zeta x.
\end{equation*}
Its discriminant $2^4\cdot3^4\cdot5^3\zeta^2$ is non-zero and $F_{2,\zeta}^+$ is an elliptic curve
over $\Fset_q$. Let us denote by $n_{2,q}^+(\zeta)$ the number of $\Fset_q$-rational points of $F_{2,\zeta}^+$ and
write
\begin{equation*}
    b_q^+(\zeta)=q+1-n_{2,q}^+(\zeta).
\end{equation*}
If $\zeta\in B^{-}(n,q)$, let us define in the same way, the cubic curve
\begin{equation*}
F_{2,\zeta}^-: y^2=x^3+\delta_{2,\zeta}x^2+5\zeta x.
\end{equation*}
As a twist of $F_{2,\zeta}^+$, it is also an elliptic curve over $\Fset_q$ and we write
\begin{equation*}
    b_q^-(\zeta)=q+1-n_{2,q}^-(\zeta),
\end{equation*}
where $n_{2,q}^-(\zeta)$ denotes the number of $\Fset_q$-rational points of $F_{2,\zeta}^-$. 
\end{enumerate}

Our criterion is stated in the following theorem which is a refinement of \cite[Th.$1.4$]{Bil07}.
\begin{thm}
Let $p$ be a prime number $\ge 17$. Assume that the following two conditions hold.
\begin{enumerate}
    \item For the curve $f=\St{1200}{K}{1}$, there exists an integer $n\ge 2$ such that
    \begin{enumerate}
        \item the integer $q=np+1$ is a prime number;
        \item we have $a_q'^2\not\equiv 4\pmod{p}$;
        \item for all $\zeta$ in $A^+(n,q)$, we have $a_q'\not\equiv a_q^+(\zeta) \pmod{p}$;
        \item for all $\zeta$ in $A^-(n,q)$, we have $a_q'\not\equiv a_q^-(\zeta) \pmod{p}$.
    \end{enumerate}
    \item For the curve $f=\St{1200}{A}{1}$, there exists an integer $n\ge 2$ such that
    \begin{enumerate}
        \item the integer $q=np+1$ is a prime number;
        \item we have $a_q'^2\not\equiv 4\pmod{p}$;
        \item for all $\zeta$ in $B^+(n,q)$, we have $a_q'\not\equiv b_q^+(\zeta)\pmod{p}$;
        \item for all $\zeta$ in $B^-(n,q)$, we have $a_q'\not\equiv b_q^-(\zeta)\pmod{p}$.
    \end{enumerate}
\end{enumerate}
Then, there is no non-trivial primitive solution of $x^5+y^5=3z^p$.
\end{thm}
\pf Let $n$ as in the theorem. By Lemma~\ref{lem:remainingcurves}, $\rho_p$ is isomorphic to the mod~$p$ representation $\overline{\sigma_{f,p}}$ of $f=\St{1200}{A}{1}$ or $\St{1200}{K}{1}$. If $E(a,b)$ does not have good reduction at $q$, then $E(a,b)$ has multiplicative reduction (\cite[Lem.$2.7$]{Bil07}) and by~\cite[prop.$3(iii)$]{KrOe92}, we have
\begin{equation*}
 a_q'\equiv \pm(q+1)\equiv \pm 2\pmod{p}. 
\end{equation*}
This contradicts the conditions~$1(b)$ and~$2(b)$ of the theorem. So, we deduce that $E(a,b)$ has good reduction at $q$ or in other words that $q$ does not divide $c$. 

We now follow step by step the discussion of~\cite[\S$4.4$]{Bil07} without giving all the details. Let us denote by $\phi$ the polynomial $\phi(x,y)=x^4-x^3y+x^2y^2-xy^3+y^4$ and by $\overline{a}$ (resp. $\overline{b}$) the reduction of $a$ (resp. $b$) modulo~$q$.
\begin{enumerate}
	\item Assume that $5$ divides $a+b$. Then, there exists $c_1$ and $c_2$ two integers such that 
\begin{equation*}
    5(a+b)=3c_1^p, \quad \phi(a,b)=5c_2^p \quad \text{and} \quad c=c_1c_2.
  \end{equation*}
Furthermore, if $u=c_1^p\pmod{q}$ and $v=c_2^p\pmod{q}$, then 
\begin{equation*}
    \overline{a}'=\frac{\overline{a}}{u},\quad \overline{b}'=\frac{\overline{b}}
{u}\quad \text{and}\quad \zeta=\frac{v}{u^4},
\end{equation*}
satisfy
\begin{equation*}
    5(\overline{a}'+\overline{b}')=3\quad\text{and}\quad \phi(\overline{a}',\overline{b}')=5\zeta.
\end{equation*}
We then deduce that $\overline{b}'$ is a root of the following polynomial
\begin{equation*}
    P_{1,\zeta}(X)=X^4-\frac{6}{5}X^3+\frac{18}{25}X^2-\frac{27}{125}X+\frac{81}{3125}-\zeta\in\Fset_q\lbrack 
    X\rbrack.
\end{equation*}
So, $\overline{b}'$ is one of the following elements
\begin{equation*}
    \frac{3}{10}+\frac{\alpha_{1,\zeta}}{50},\quad\frac{3}{10}-\frac{\alpha_{1,\zeta}}{50},\quad
    \frac{3}{10}+\frac{\beta_{1,\zeta}}{50},\quad\frac{3}{10}-\frac{\beta_{1,\zeta}}{50},
\end{equation*}
where $\alpha_{1,\zeta}$ (resp. $\beta_{1,\zeta}$) is a square root of $-225+10\delta_{1,\zeta}$ (resp. $-225-10\delta_{1,\zeta}$) modulo~$q$.

\begin{enumerate}
	\item Assume that we have 
\begin{equation*}
\left\{\overline{a}',\overline{b}'\right\}=\left\{\frac{3}{10}+\frac{\alpha_{1,\zeta}}{50},\frac{3}{10}-\frac{\alpha_{1,\zeta}}{50}\right\}.
\end{equation*}
Then, $\zeta$ belongs to the set $A^+(n,q)$ and the reduction modulo~$q$ of the curve $E(a,b)$ is isomorphic to $F_{1,\zeta}^+$. So, we deduce that 
\begin{equation*}
a_q \equiv a_q^+(\zeta) \pmod{p}.
\end{equation*}
But, by Lemma~\ref{lem:remainingcurves}, we have $a_q\equiv a_q'\pmod{p}$, where $a_q'$ is the $q$th Fourier coefficient of $\St{1200}{K}{1}$. This contradicts our hypothesis~$1(c)$. 
	\item Assume that we have 
\begin{equation*}
\left\{\overline{a}',\overline{b}'\right\}=\left\{\frac{3}{10}+\frac{\beta_{1,\zeta}}{50},\frac{3}{10}-\frac{\beta_{1,\zeta}}{50}\right\}.
\end{equation*}
Then, $\zeta$ belongs to the set $A^-(n,q)$ and the reduction modulo~$q$ of the curve $E(a,b)$ is isomorphic to $F_{1,\zeta}^-$. So, we deduce that 
\begin{equation*}
a_q \equiv a_q^-(\zeta) \pmod{p}.
\end{equation*}
But, by Lemma~\ref{lem:remainingcurves}, we have $a_q\equiv a_q'\pmod{p}$, where $a_q'$ is the $q$th Fourier coefficient of $\St{1200}{K}{1}$. This contradicts our hypothesis~$1(d)$.
\end{enumerate}
We finally deduce that $5$ does not divide $a+b$.

	\item If $5$ does not divide $a+b$, then there exists $c_1$ and $c_2$ two integers such that 
\begin{equation*}
    a+b=3c_1^p, \quad \phi(a,b)=c_2^p \quad \text{and} \quad c=c_1c_2.
  \end{equation*}
Furthermore, if $u=c_1^p\pmod{q}$ and $v=c_2^p\pmod{q}$, then 
\begin{equation*}
    \overline{a}'=\frac{\overline{a}}{u},\quad \overline{b}'=\frac{\overline{b}}
{u}\quad \text{and}\quad \zeta=\frac{v}{u^4},
\end{equation*}
satisfy
\begin{equation*}
    \overline{a}'+\overline{b}'=3\quad\text{and}\quad \phi(\overline{a}',\overline{b}')=\zeta.
\end{equation*}
We then deduce that $\overline{b}'$ is a root of the following polynomial
\begin{equation*}
   P_{2,\zeta}(X)=X^4-6X^3+18X^2-27X+\frac{81-\zeta}{5}\in\Fset_q\lbrack X\rbrack.
\end{equation*}
So, $\overline{b}'$ is one of the following elements
\begin{equation*}
    \frac{3}{2}+\frac{\alpha_{2,\zeta}}{10},\quad\frac{3}{2}-\frac{\alpha_{2,\zeta}}{10},\quad
    \frac{3}{2}+\frac{\beta_{2,\zeta}}{10},\quad\frac{3}{2}-\frac{\beta_{2,\zeta}}{10},
\end{equation*}
where $\alpha_{2,\zeta}$ (resp. $\beta_{2,\zeta}$) is a square root of $-225+10\delta_{2,\zeta}$ (resp. $-225-10\delta_{2,\zeta}$) modulo~$q$.

\begin{enumerate}
	\item Assume that we have 
\begin{equation*}
\left\{\overline{a}',\overline{b}'\right\}=\left\{\frac{3}{2}+\frac{\alpha_{2,\zeta}}{10},\frac{3}{2}-\frac{\alpha_{2,\zeta}}{10},\right\}.
\end{equation*}
Then, $\zeta$ belongs to the set $B^+(n,q)$ and the reduction modulo~$q$ of the curve $E(a,b)$ is isomorphic to $F_{2,\zeta}^+$. So, we deduce that 
\begin{equation*}
a_q \equiv b_q^+(\zeta) \pmod{p}.
\end{equation*}
But, by Lemma~\ref{lem:remainingcurves}, we have $a_q\equiv a_q'\pmod{p}$, where $a_q'$ is the $q$th Fourier coefficient of $\St{1200}{A}{1}$. This contradicts our hypothesis~$2(c)$. 
	\item Assume that we have 
\begin{equation*}
\left\{\overline{a}',\overline{b}'\right\}=\left\{\frac{3}{2}+\frac{\beta_{2,\zeta}}{10},\frac{3}{2}-\frac{\beta_{2,\zeta}}{10}\right\}.
\end{equation*}
Then, $\zeta$ belongs to the set $B^-(n,q)$ and the reduction modulo~$q$ of the curve $E(a,b)$ is isomorphic to $F_{2,\zeta}^-$. So, we deduce that 
\begin{equation*}
a_q \equiv b_q^-(\zeta) \pmod{p}.
\end{equation*}
But, by Lemma~\ref{lem:remainingcurves}, we have $a_q\equiv a_q'\pmod{p}$, where $a_q'$ is the $q$th Fourier coefficient of $\St{1200}{A}{1}$. This contradicts our hypothesis~$2(d)$.
\end{enumerate}
We finally deduce that there is no non-trivial primitive solution of the equation $x^5+y^5=dz^p$.
\end{enumerate}

\medskip

\rmk For a given $p$, a \verb"pari/gp" program giving an integer $n$ as in the theorem is available at the address \newline\verb"http://www.institut.math.jussieu.fr/~billerey/Fermatnew".

\end{document}